\documentclass[11pt]{amsart}

\usepackage[latin1]{inputenc}
\usepackage[T1]{fontenc}
\usepackage{lmodern}

\usepackage{amssymb}

\newcommand{\R}{\mathbb{R}}
\newcommand{\Z}{\mathbb{Z}}
\newcommand{\Q}{\mathbb{Q}}

\newcommand{\g}{{\mathfrak g}}
\newcommand{\F}{\mathcal{F}}
\newcommand{\inv}{^{-1}}

\newcommand{\id}{\mathrm{id}}
\newcommand{\image}{\mathop{\mathrm{im}}}

\newcommand{\DIFF}[1]{ \mathop{\mathrm{ Diff}}(#1)}
\newcommand{\AUT}[2]{{\mathrm{Aut}_{#2}(#1)}}

\newcommand{\DIFFGA}{{\mathrm{Diff}}(G : \Gamma)}

\newtheorem{theo}{Theorem}[section]
\newtheorem{defi}[theo]{Definition}
\newtheorem{prop}[theo]{Proposition}
\newtheorem{exem}[theo]{Example}
\newtheorem{lemm}[theo]{Lemma}
\newtheorem{coro}[theo]{Corollary}

\theoremstyle{remark}
\newtheorem{rema}[theo]{Remark}

\title
[The diffeomorphism group of a Lie foliation]
{The diffeomorphism group of a Lie foliation}

\author{G.~Hector}
\address{
Institut Camille Jordan\\
Universit\'e Claude Bernard Lyon 1\\
69622-Villeurbanne (France)
}
\email{gilb.hector@gmail.com}

\author{E.~Mac\'{\i}as-Virg\'os}
\address{Dpto. Xeometria e Topoloxia \\
Facultade de Matem\'aticas\\
Universidade de Santiago de Compostela\\
15782-Santiago de Compostela (Spain)
}
\email{quique.macias@usc.es}
\thanks{Second author partially supported by FEDER and Research Project MTM2008-05861 MICINN Spain}

\author{A.~Sotelo-Armesto}
\address{Dpto. Xeometria e Topoloxia\\
Facultade de Matem\'aticas \\
Universidade de Santiago de Compostela\\
15782-Santiago de Compostela (Spain)
}


\begin{document}

\begin{abstract}
 We explicitly compute the diffeomorphism group of several types
of linear foliations (with dense leaves)  on the torus $T^n$, $n\geq 2$, namely codimension one foliations,
flows, and the so-called non-quadratic foliations. We show in particular that  non-quadratic foliations are rigid, in the sense that they do not admit
transverse diffeomorphisms other than $\pm \id$ and translations. The computation is an application of a general formula that we prove for the
diffeomorphism group of any Lie foliation with dense leaves on a compact manifold. Our results  generalize those  of P.~Donato and P.~Iglesias for $T^2$, P.~Iglesias and G.~Lachaud 
  for codimension one foliations on $T^n$, $n\geq 2$, and  B.~Herrera  for transcendent foliations.
The theoretical setting of the paper is that of J.~M.~Souriau's diffeological spaces. 
\end{abstract}

\maketitle

\section{Introduction}
\label{intro}
J. M. Souriau's diffeological spaces \cite{SOURIAU} are a fruitful generalization of
manifolds. In that setting one has natural constructions for
subspaces, quotient spaces and spaces of maps. In particular,
the factor space $G/\Gamma$ of the simply connected Lie group $G$ by any  
totally disconnected subgroup $\Gamma$
can be endowed with a diffeology,
and its group of diffeomorphisms $\DIFF{G/\Gamma}$  is a diffeological group ---the
generalization of the notion of Lie group.
When $\Gamma$ is dense in $G$, we prove that the group $\DIFF{G/\Gamma}$ is isomorphic to $({\AUT{G}{\Gamma}} \ltimes
G)/\Gamma$. This formula will be our main tool.

For a $G$-Lie foliation  \cite{FEDIDA-2,FEDIDA,MOLINO} with holonomy $\Gamma$ on the compact manifold $M$,
the space of leaves $M/\mathcal{F}$ turns out to be diffeomorphic (as a diffeological space) to $G/\Gamma$, hence the group
$\DIFF{M/\mathcal{F}}$ of transverse diffeomorphisms can be computed by means of the formula above (when the leaves are
dense). An important particular case (with
$G=\R^n$) is that of linear foliations on the torus $T^n$.  

In this paper we explicitly compute the  diffeomorphism group $\DIFF{M/\mathcal{F}}$ for several types of linear foliations
(with dense leaves)  on the torus
$T^n$, $n\geq 2$, namely codimension one foliations, flows, and the so-called non-quadratic foliations. We show in particular that  non-quadratic foliations are rigid, in the sense
that they do not admit transverse diffeomorphisms other than $\pm \id$ and translations. Our results
generalize those of P.~Donato and P.~Iglesias \cite{DONATOIGLESIAS} for $T^2$, P.~Iglesias and G.~Lachaud \cite{IGLESIASLACHAUD} for codimension one foliations on $T^n$, $n\geq 2$, and  B.~Herrera \cite{HERRERA-2,HERRERA} for transcendent foliations. 

\section{Diffeological spaces}
The notion of {\it diffeological space} is due to J. M. Souriau \cite{SOURIAU}. We shall briefly review the basic constructions on the category of diffeological spaces that we need along
this paper. Proofs are left to the reader, see also \cite{HECTORMACIAS,HECTOR,IGLESIASBOOK}.

Let $M$ be a set. Any set map $\alpha \colon U\subset \R^n \to M$ defined
on an open set $U$ of some $\R^n, n\geq 0,$ will be called a {\it plot}
on $M$. The name {\it plot} is chosen instead of
{\it chart} to avoid some confusion with the usual notion of chart in a manifold.  When possible, a plot
$\alpha
$ with domain
$U$ will be simply denoted by $\alpha_U$.

\begin{defi} \label{DIFEOL}
A {\em diffeology} of class $\mathcal{C}^\infty$ on the set $M$ is any collection
$\mathcal{P}$ of plots $\alpha  \colon U_\alpha \subset \R^{n_\alpha}\to M$,
$n_\alpha \geq 0$, verifying the following axioms:
\begin{enumerate}
\item
Any constant map $c \colon  \R^n \to M$, $n\geq 0$, belongs to
$\mathcal{P}$;
\item
Let $\alpha \in \mathcal{P}$ be defined on $U
\subset \R^n$ and let $h \colon V \subset \R^m \to U
\subset \R^n$ be any $\mathcal{C}^\infty$ map; then $\alpha \circ h \in \mathcal{P}$;
\item
Let $\alpha \colon U \subset \R^ n \to M$ be a plot. If  any
$t\in U$ has a neighbourhood $U_t$ such that $\alpha_{\mid U_t}$ belongs to $\mathcal{P}$ then
$\alpha \in \mathcal{P}$.
\end{enumerate}
\end{defi} 

Usually, a diffeology $\mathcal{P}$ on the set $M$ is defined by means of a {\it generating set}, that is
by giving any set $\mathcal{G}$  of plots (which is implicitly supposed to contain all constant maps)
and taking the least diffeology containing it. Explicitly,  the diffeology
$\left<\mathcal{G}\right>$ generated by $\mathcal{G}$ is the set of plots $\alpha \colon U \to M$
such that any point $t\in U$ has a neighbourhood $U_t$ where $\alpha$ can be written as $\gamma\circ h$
for some $\mathcal{C}^\infty$ map $h$ and some $\gamma\in \mathcal{G}$.

\begin{exem}  A finite dimensional manifold $M$  is endowed with the diffeology generated by the
charts
$U\subset \R^n \to M$, $n=\dim M$, of any atlas.
\end{exem}

\subsection{Basic constructions} \label{BASICS}

A map $F\colon (M,\mathcal{P}) \to (N,\mathcal{Q})$ between diffeological spaces is {\em differentiable}
if $F\circ \alpha \in \mathcal{Q}$ for all $\alpha \in \mathcal{P}$. A {\em diffeomorphism} is a differentiable
map with a differentiable inverse.

Let $(M,\mathcal{P})$ be a diffeological space  and $F\colon M \to N$ a map of sets. The {\em final
diffeology} $F_\star\mathcal{P}$ on $N$  is that generated by the plots $F\circ\alpha$, $\alpha\in\mathcal{P}$.
A particular case is the {\em quotient diffeology} associated to an equivalence relation on $M$.

Analogously, let $(N,\mathcal{Q})$ be a diffeological space and $F\colon M \to N$ a map of sets. The
{\em initial diffeology} $F^\star\mathcal{Q}$ on $M$ is that generated by the plots $\alpha$ in $M$
such that $F\circ \alpha \in \mathcal{Q}$. A particular case is  the {\em induced diffeology}
on any subset $M\subset N$.

Both constructions verify the usual universal properties.

Let $(M,\mathcal{P})$, $(N,\mathcal{Q})$ be two diffeological spaces. We can endow the cartesian product $M\times N$
with the {\it product diffeology} $\mathcal{P}\times\mathcal{Q}$ generated by the plots $\alpha\times\beta$, $\alpha\in\mathcal{P}$,
$\beta\in\mathcal{Q}$.

Finally, let $\mathcal{D}(M,N)$ be the space of differentiable maps between two diffeological spaces $(M,\mathcal{P})$
and $(N,\mathcal{Q})$. We define the {\em functional diffeology} on it by taking as a generating
set all plots $\alpha \colon U \to  \mathcal{D}(M,N)$ such that the associated map
$\widehat{\alpha} \colon U\times M \to N$ given by $\widehat{\alpha}(t,x) =\alpha(t)(x)$ is
differentiable.  

\subsection{Diffeological groups}
\begin{defi}  A diffeological group is a diffeological space $(G,\mathcal{P})$ endowed with a group structure
such that the division map $\delta \colon G\times G \to G$, $\delta(x,y)= x y\inv$, is differentiable.
\end{defi} 
A typical example of diffeological group is the diffeomorphism group of a finite dimensional manifold $M$, endowed with the diffeology induced
by $\mathcal{D}(M,M)$ \cite{HECTOR}.
In
\cite{HECTORMACIAS} (see also Prop. \ref{DIFF}) it is proven that the diffeomorphism group of the space of leaves of a Lie
foliation   is a diffeological group too. The aim of the present paper is to compute it.

\section{Homogeneous spaces and covering maps}
Let $G$ be a {\it connected} and {\it simply connected} Lie group, $\Gamma\subset G$ an arbitrary {\em totally disconnected} subgroup.  While it is possible to develop a general theory of fibre bundles and covering spaces in the diffeological
setting \cite{IGLESIAS}, we shall directly prove some lifting properties for the quotient map $\pi \colon G \to G/\Gamma$.
These lifting properties  will be essential later. 
\subsection{Lifting of diffeomorphisms} \label{LIFTING}

The factor space $G/\Gamma$ is endowed with the quotient
diffeology, that is the collection of plots $\alpha \colon U \to G/\Gamma$ which locally lift through $\pi$ to a smooth
map $U \to G$.

\begin{prop}  \label{LIFT}  Any differentiable map (in the diffeological sense)
$\overline\varphi
\colon {G/\Gamma } \to {G/\Gamma }$ has a smooth
lift $\varphi \colon G \to G$, that is a $\mathcal{C}^\infty$ map such that
$\pi\varphi = \overline\varphi \pi$.
\end{prop} 

\proof  By definition of quotient diffeology, the result is {\it locally} true, that is for any $x\in G$
there exists an open neighbourhood $U_x$ and a smooth map $\varphi _x \colon U_x \to G$
such that $\pi \circ \varphi _x = \overline\varphi \circ \pi$ on $U_x$. We can suppose that $U_x$ is a connected open set.

Now we  define an integrable distribution $\mathcal{D}$ on $G \times G$ in the following way. Since an arbitrary point $(x,g)\in G \times G$  can be
written  as $(x,\varphi _x(x)h)$ for some $h\in G$, let
$$\mathcal{D}_{(x,g)}=\{(v,(R_h\circ \varphi_x)_{*x}(v) )\colon v\in T_xG\} \subset T_{(x,g)}(G\times G).$$
The distribution $\mathcal{D}$ is well defined because  two local lifts differ by some translation. 
In fact, let $x,y\in G$ such that $U_x \cap U_y\neq\emptyset$. Then for any $z\in U_x \cap U_y$ and any
connected neighbourhood $V_z \subset U_x \cap U_y$, the local lifts $\varphi _x, \varphi _y$  define the  continuous map $\gamma \colon V_z \to \Gamma$ given by $\gamma(t) =\varphi_x(t)\inv\varphi_y(t)$. Since
 $V_z$ is connected
and  the set $\Gamma$ is totally disconnected, the map $\gamma$ must be constant, hence $\varphi_y=R_\gamma\circ\varphi_x$.

Moreover $\mathcal{D}$ has constant rank, and it is integrable, the integral submanifolds being translations of the graphs of the local lifts.

 Let us choose   some point
$x_0\in G$ such that $[x_0] =  \overline\varphi([e])$, and let $\widetilde{G}$ be the  maximal integral
submanifold passing through $(e,x_0)$. We shall prove in the next Lemma \ref{CUB} that the projection of
$\widetilde G \subset G\times G$  onto the first factor is a covering map. Since $G$ is simply connected, it follows that $\widetilde G$ is the graph of a global lift. \qed

\begin{lemm}  \label{CUB} The projection $p_1\colon \widetilde G \subset G\times G \to G$ is a covering map.
\end{lemm} 

\proof Clearly $p_1$ is a differentiable submersion, hence an open map,  so $p_1(\widetilde G)$ is an open
subspace of
$G$. Let us   prove that it is closed too; this will show that the map $p_1\colon \widetilde G \to G$
is onto, because the Lie group $G$ is connected. 

Suppose $x\in G$ is   in the closure of
$p_1(\widetilde G)$, and let $U_x$ be a connected open neighbourhood where  the local lift $\varphi_x$  is defined. Let $y\in U_x \cap p_1(\widetilde G)$, then $(y,\varphi_x(y)h) \in \widetilde G$ for some $h\in G$. This   implies that the graph of $R_h\circ \varphi_x$, which is an integral submanifold of $\mathcal{D}$,  is contained in $\widetilde G$. Hence $x\in p_1(\widetilde G)$.  

It remains to prove  that any $x\in G$ has a neighbourhood $U_x$ such that $(p_1)\inv(U_x)$ is
a disjoint union of open sets, each one homeomorphic to $U_x$ by $p_1$. It is clear that we
can restrict ourselves to the case $x=e$. Let $\varphi _e \colon U_e \to G$ be a connected local lift of
$\overline\varphi$. We can suppose that $\varphi _e(e) =x_0$. Let $\widetilde{U_e}$ be its graph. Then $\widetilde{U_e}$ is an open
subset of $\widetilde G$, containing $(e,x_0)$, with $p_1(\widetilde{U_e})=U_e$. Let $I$ be the non-empty set
$$I=\{\gamma \in \Gamma \colon (e,x_0\gamma) \in \widetilde G\}.$$
Then $(p_1)\inv(U_e)$ is
the disjoint union of  the sets $R_\gamma(\widetilde{U_e})$,  $\gamma\in I$.  \qed

\begin{coro}  Any diffeomorphism  of $G/\Gamma$ can be lifted to a diffeomorphism  of
$G$.
\end{coro} 

\proof In the proof of Theorem \ref{LIFT}, the integral manifold $\widetilde G$ passing through a given point is
unique, hence the lift of the map $\overline\varphi$ is unique once one has selected a representative of
$\overline\varphi([e])$. Then the appropriate lifts of $\overline\varphi$ and $(\overline\varphi)\inv$ are
inverse maps. \qed

Analogous arguments give the following result.

\begin{coro} \label{LIFTOPEN2} Let $U$ be a connected simply connected open subset of $\R^n$,
$n\geq 0$. Any differentiable map (resp. diffeomorphism) $  U\times {G/\Gamma } \to U\times {G/\Gamma }$ can
be lifted to a $\mathcal{C}^\infty$ map (resp. diffeomorphism) $ U\times G \to U\times G$.
 \end{coro}

 \section{The diffeomorphism group}
Several results in this section were previously announced by two of the authors in   \cite{HECTORMACIAS}.
For the sake of completeness,  we shall give
complete proofs.

Let   ${G/\Gamma } $ be the factor space (endowed with
the quotient diffeology) of the connected simply connected Lie group $G$ by the totally disconnected subgroup $\Gamma$.  
Let $\DIFF{{G/\Gamma } }$ be the diffeomorphism group of $G/\Gamma$, with
the diffeology induced by the functional diffeology of $\mathcal{D}({G/\Gamma } , {G/\Gamma } )$.

\begin{lemm} \label{LIFTOPEN1}
Any plot  in
$\DIFF{{G/\Gamma } }$ with simply connected domain can be lifted to some plot in $\DIFF{G}$.
 \end{lemm} 

\proof Let
$\alpha \colon U \subset \R^n \to
\DIFF{{G/\Gamma } } \subset \mathcal{D}({G/\Gamma }, {G/\Gamma } )$ be a
plot in the functional diffeology, that is such that the map
$$\widehat{F}\colon U\times {G/\Gamma } \to U \times {G/\Gamma },\quad
\widehat{F}(t,[g]) = (t,\alpha(t)([g])),$$ is differentiable. Take (Corollary \ref{LIFTOPEN2}) a lift $F\colon
U\times G
\to U\times G$. Then    $F(t,g) = (t, \xi(t,g))$, where
$[\xi(t,g)] = \alpha(t)([g])$. For each $t\in U$,  lift $\alpha(t)\inv$  to
some diffeomorphism $\lambda_t$ of the Lie group $G$. Since $[(\lambda_t\circ \xi)(t,g)] =[g]$ for
all $g\in G$, we have $\lambda_t \circ \xi(t,-) = R_\gamma$ for some $\gamma \in \Gamma$,
because the Lie group $G$ is connected  and the subgroup $\Gamma$ is  totally disconnected.
Analogously
$\xi(t,-)
\circ \lambda_t =R_\mu$ for some $\mu \in \Gamma$. This proves that $ \xi(t,-)$ is a
diffeomorphism of $G$, so $\xi \colon U \to \DIFF{G}$ is the desired lift.
\qed

\begin{theo}  [\cite{HECTORMACIAS}]\label{DIFF}
$\DIFF{{G/\Gamma } }$ is a diffeological group, with the diffeology induced by $\mathcal{D}(G/\Gamma, G/\Gamma)$.
\end{theo}  

\proof The differentiability of the composition follows from the very definitions (for the diffeomorphism
group of any diffeological space). The crucial point is to prove that the inversion map $I\colon \DIFF{{G/\Gamma } } \to \DIFF{{G/\Gamma } }$
is differentiable.

Let $\alpha_U$ be a plot in $\DIFF{{G/\Gamma } }$. We can suppose that the domain $U$ is simply connected. By applying the
inverse function theorem to the manifold
$U\times G$ we conclude that the map $F$ in Lemma \ref{LIFTOPEN1} is a diffeomorphism.
But   $F\inv$ is a lifting of the map $(t,g)\in U\times {G/\Gamma } \mapsto
(t,\alpha(t)\inv([g])) \in U\times {G/\Gamma }$, which proves that $I\circ\alpha$ is a 
plot in the functional diffeology. \qed

\subsection{An explicit formula}\label{FORMULA}
Let us suppose that the subgroup $\Gamma$ is  {\it dense} in $G$.

We denote ${\AUT{G}{\Gamma}} \subset \AUT{G}{}$  the group of automorphisms of the Lie group $G$
which preserve the subgroup $\Gamma$. Let us consider the semidirect product ${\AUT{G}{\Gamma}} \ltimes G$, where the product is  given by
$(\theta_1,g_1)(\theta_2,g_2) = (\theta_1\theta_2, g_1\theta_1(g_2))$. The group $\Gamma$
can be identified with a subgroup of that semidirect product by means of the map
$i(\gamma) = (i_\gamma, \gamma\inv)$, where $i_\gamma$ is the inner automorphism
$i_\gamma(g) = \gamma g \gamma\inv$.  This subgroup is {\it invariant}, because 
$$(\theta,g)(i_\gamma,\gamma\inv)(\theta,g)\inv = (i_{\theta(\gamma)},\theta(\gamma)\inv).$$

\begin{theo}   \label{MAIN} Let $G$ be a connected simply connected
Lie group, $\Gamma \subset G$ a totally disconnected dense subgroup. Then $\DIFF{G/\Gamma}$ is isomorphic to the
quotient group
$({\AUT{G}{\Gamma}} \ltimes G)/\Gamma$.
\end{theo}  

This formula  will be our main tool.
The proof (sketched in \cite{HECTORMACIAS}) will be divided into several Lemmas.

\begin{lemm}  \label{MOR} Let $\varphi$ be any diffeomorphism of $G$ which induces a
diffeomorphism on ${G/\Gamma }$. Let us define
$\theta_\varphi = L_{\varphi(e)}\inv \circ
\varphi$. Then $\theta_\varphi \in {\AUT{G}{\Gamma}}$. Moreover $\theta_{\varphi\circ\psi}=\theta_{\varphi}\circ\theta_{\psi}$.
\end{lemm} 

\proof We have $\varphi(g) = \varphi(e)\theta_\varphi(g)$, for all $g\in G$. But in fact $\varphi(g\gamma) =
\varphi(g)\theta_\varphi(\gamma)$ for all $\gamma\in\Gamma$, because for any fixed $\gamma$ the
continuous map $\varphi(g)\inv\varphi(g\gamma)$ defined on $G$ takes its values in the
totally disconnected group $\Gamma$, hence it is constant. Moreover
$\theta_\varphi(\gamma)\in\Gamma$.

Now, for $\gamma,\mu \in \Gamma$ we have
$$\theta_\varphi(\mu\gamma) =
\varphi(e)\inv\varphi(\mu\gamma) =\varphi(e)\inv\varphi(\mu)\theta_\varphi(\gamma) =
\theta_\varphi(\mu)\theta_\varphi(\gamma),$$
 hence $\theta_\varphi$ is a group morphism when restricted to $\Gamma$.
Since
the subgroup $\Gamma \subset G$ is  dense,   it follows that $\theta_\varphi$ is an automorphism
of
$G$.

Finally, for two diffeomorphisms $\varphi,\psi$ we have
$$\varphi\psi = L_{\varphi(e)}\theta_\varphi L_{\psi(e)}\theta_\psi =
 L_{\varphi(e)}L_{\theta_\varphi\psi(e)}\theta_\varphi\theta_\psi =
L_{\varphi\psi(e)}\theta_\varphi\theta_\psi$$ which proves
$\theta_{\varphi\psi}=\theta_{\varphi}\theta_{\psi}$. 
\qed

\begin{rema}
Clearly, $\theta_\varphi =\id$ if and only if $\varphi$ is a (left) translation by an element
of
$G$.  On the other hand, right translations by elements of $\Gamma$ induce the identity on
${G/\Gamma }$.
\end{rema} 

We know from Theorem \ref{LIFT} that any diffeomorphism $\overline\varphi$ of $G/\Gamma$ can be lifted
to some diffeomorphism $\varphi$ of $G$.

\begin{lemm} 
The map $\Phi \colon \DIFF{G/\Gamma} \to  ({\AUT{G}{\Gamma}} \ltimes G)/\Gamma$ given by
$$\Phi(\overline\varphi) = [(\theta_\varphi, \varphi(e))]$$
is an isomorphism of groups.
\end{lemm} 

\proof That the map $\Phi$ is well defined follows from the fact that if $\varphi, \psi$
are two different lifts
 of
$\overline\varphi$, then
 $\theta_\psi = i_\gamma \circ \theta_\varphi$, where $i_\gamma$ is the inner automorphism with $\gamma =
\psi(e)\inv\varphi(e)$, hence
we have
$$(\theta_\psi,\psi(e)) = (i_\gamma,\gamma\inv)(\theta_\varphi,\varphi(e)).$$

Let us prove that $\Phi$ is a bijective map. First it is injective, because if $\overline\varphi$
goes
into the identity, then $\theta_\varphi = i_\gamma$ and $\varphi(e) = \gamma\inv$, for some
$\gamma\in\Gamma$, which implies $\varphi  = R_\gamma\inv$,
that is
$\overline\varphi$ is the identity  of $G/\Gamma$.  On the other hand, $\Phi$ is onto because for
any pair $(\theta,g) \in {\AUT{G}{\Gamma}} \ltimes G$ we can take the diffeomorphism $\varphi = L_g\circ \theta$.

That $\Phi$ is a morphism of groups follows from  Lemma \ref{MOR}. \qed

\begin{lemm}  $\Phi$ is a diffeomorphism of diffeological spaces.
\end{lemm} 

\proof From the definition of induced diffeology (Subsection \ref{BASICS}) it turns out that any
subgroup of a diffeological group is a diffeological group. This applies to ${\AUT{G}{\Gamma}} \subset
\DIFF{G}$, and also to the subgroup $\DIFFGA$ of diffeomorphisms of $G$ which induce a
diffeomorphism
on ${G/\Gamma }$. We have the following diagram of exact sequences of diffeological groups:
$$\begin{matrix}
Z(\Gamma) & \to & G & = & G/Z(\Gamma) \cr
\downarrow && L \downarrow && L \downarrow \cr
\Gamma & \stackrel{R}{\to} & \DIFFGA & \to & \DIFF{{G/\Gamma }} \cr
i \downarrow & &\theta \downarrow && \downarrow  \cr
i(\Gamma) & \to & {\AUT{G}{\Gamma}} & \to & {\AUT{G}{\Gamma}} / i(\Gamma)\cr
\end{matrix}
$$
where $i(\gamma)$ is an inner automorphism, $R(\gamma)$ is a right translation, and $L(g)$ is a
left translation. Then the differentiability of $\Phi$ and $\Phi\inv$ is an easy consequence
of Lemma \ref{LIFTOPEN1} and the diagram above. \qed

This completes the proof of Theorem \ref{MAIN}.

\section{Lie foliations}\label{LIE}
Lie foliations play a central role in the study of transversely riemannian foliations \cite{MOLINO}.
Let $\g$ be a Lie algebra of dimension $n$. A $\g$-Lie foliation on the  manifold $M$ is defined
as the kernel of a non-degenerate $1$-form $\omega$ with values in $\g$,  such that $d\omega= (-1/2)[\omega,\omega]$.
Once a basis $X_1,\dots,X_n$ of $\g$ with structural constants $c_{ij}^k$ has been fixed,  this is equivalent to having
$n$ independent  real
$1$-forms $\omega_1,\dots,\omega_n$ on $M$ such that $d\omega_k = \sum{c_{ij}^k\omega_i\wedge\omega_j}$.

When the manifold $M$ is {\it compact}, it is  well known  \cite{FEDIDA-2,MACIAS}
that there exists a regular covering $p\colon
\widetilde M \to M$ such that the lifted foliation $p^\star\mathcal{F}$ is  a locally trivial bundle
$D \colon\widetilde M \to G$ over the connected simply connected Lie group $G$ associated to $\g$. By fixing
base points $x_0\in M$ and $\widetilde x_0\in\widetilde M$ with $p(\widetilde x_0)=x_0$ and
$D(\widetilde x_0) = e\in G$, one obtains a group morphism $h \colon\pi_1(M) \to G$ such
that  $\ker h = \pi_1(\widetilde M)$. Moreover $D$ is $h$-equivariant. The {\it developing map} $D$ and the {\it
holonomy morphism}
$h$ completely determine the foliation. The image $\Gamma$ of $h$ is called the {\it holonomy group}; it is dense in $G$ if and only if all the leaves of $\mathcal{F}$ are dense in $M$.

From the definitions above it is easy to prove that the space of leaves $M/\mathcal{F}$ (endowed
with the quotient diffeology) is diffeomorphic to the factor space $G/\Gamma$.

\begin{prop}  For any Lie foliation $\mathcal{F}$ on the manifold $M$, the groups $\DIFF{M/\mathcal{F}}$  and  $\DIFF{{G/\Gamma }}$ are isomorphic.
\end{prop} 
When $M$ is compact, the group $\Gamma=\image{h}$ is finitely generated, so we can apply our previous computation of the diffeomorphism group.

\subsection{Linear foliations on the torus }
A particular case of Lie foliation is that of a linear foliation on the torus. An arbitrary linear foliation $\mathcal{F}$ of dimension
$m\geq 1$ and codimension $n\geq 1$ on the torus $M=T^{m+n}$  is determined by  some
linear subspace
$V\subset \R^{m+n}$ of dimension $m$. The holonomy group $\Gamma$ is a finitely generated subgroup of $\R^n$, hence a free abelian group of rank $k$. 
Let $V^\perp$ be the orthogonal subspace to $V$, and take arbitrary basis $v_1,\dots,v_m$ of $V$ and $w_1,\dots,w_n$ of $V^\perp$.
Then the closed  $1$-form $\omega
=(dw_1,\dots,dw_n)$ with values in $\R^n$  defines the foliation  $\mathcal{F}$. 

The holonomy morphism $\pi_1(T^{m+n})={\Z}^{m+n} \to {\R}^n$ is obtained by computing the group of periods of the form $\omega$, then $\Gamma =D({\Z}^{m+n})$ if we take  as a developping map  the orthogonal projection $D\colon {\R}^{m+n} \to V^\perp = {\R}^n$. Note that in general $V\cap {\Z} ^{m+n}\neq 0$, which means that the covering $\widetilde M$ of Section \ref{LIE} is not the universal covering $\R^{m+n}$ but the intermediate covering  corresponding to $\ker h$.

Our general formula in \ref{FORMULA} reduces to  $\DIFF{T^{m+n}/{\F}} = 
\AUT{\R^n}{\Gamma}\times (\R^n/\Gamma)$  because
the group $G=\R^n$ is abelian.

\begin{rema}
In \cite{AZIZ, AZIZ-2} A.~El Kacimi and A.~Tihami  computed the bigraded cohomology of linear
foliations on tori.
\end{rema} 

\subsection{Duality}  We can establish some kind of duality between the foliations defined by $V$ and $V^\perp$.
This idea will play an important role in our paper.
First, in order to compute the  group of transverse diffeomorphisms we shall  suppose that the   leaves of the foliation are simply connected (planes). If this is not the case,  we just need to reduce the dimension of the ambient torus. This condition   is equivalent to $V\cap \Z^{m+n}=0$, which means that  the rank $k$ of $\Gamma$ equals $m+n$, because $D\colon \Z^{m+n}\subset \R^{m+n} \to \Gamma \subset \R ^n$ is an isomorphism.   

\begin{lemm}  \label{CURIOSO} A finitely generated subgroup $\Gamma$ of $\R^p$ is dense if and only if $f(\Gamma)$ is dense
in $\R$ for any non-degenerate linear form $f \colon \R^p \to \R$.
\end{lemm} 

We leave the proof to the reader.

\begin{prop}  $V\cap \Z^{m+n}=0$ if and only if the orthogonal foliation $\F^\perp$ determined by $V^\perp$ has dense leaves.
\end{prop} 

\proof 
Let $D^\perp \colon \R^{m+n} \to V=\R^m$ be the orthogonal projection onto $V$. The holonomy group of $\F^\perp$ is $\Gamma^*=D^\perp(\Z^{m+n})$.  If $\Gamma^*$ is not dense in $\R ^m$, then by Lemma \ref{CURIOSO} there exists some non trivial linear map $f\colon \R^{m} \to \R$ such that $f(\Gamma^*)$ is not dense in $\R$, hence it is a non-zero  discrete subgroup. By multiplying by some scalar we can suppose that $f(\Gamma^*)=\Z$. Composing $f$ with $D^\perp$ then gives a linear map $\varphi \colon \R ^{m+n} \to \R$ such that $\varphi(V ^\perp)=0$ and $\varphi(\Z ^{m+n})=\Z$. Let $v\in V$, $v\neq 0$,  such that $\varphi =
\left\langle {v,-} \right\rangle$. Then $\left\langle {v,\Z^{m+n}} \right\rangle \subset \Z$ means that $v$ has integer coordinates with respect to the canonical basis of $\R ^{m+n}$, that is $v\in V\cap \Z^{m+n}\neq 0$.

The reciprocal is completely analogous.\qed
 
\begin{coro}  Suppose that $\F$ has simply connected dense leaves. Then so has $\F^\perp$, and
the groups $\AUT{\R^n}{\Gamma}$ and $\AUT{\R^m}{\Gamma^*}$ are isomorphic.
\end{coro} 

\proof 
Computing $\AUT{\R^n}{\Gamma}$  is equivalent to finding the matrices $A\in PSL(m+n,\Z)$ such that $A(V)=V$. Clearly
 this is equivalent to the condition $A^T(V^\perp)=V^\perp$ when  $V\cap \Z^{m+n}=0$ and $V^\perp \cap \Z^{m+n}=0$.
\qed

Later in Section \ref{FLOWS} we shall use  the fact that this duality always occurs when $\mathcal{F}$ is a dense flow.

\begin{rema}
All along the following sections we shall consider the  following {\it Moebius action} of $GL(k,\Z)$ on the
projective space
$P\R^{k-1}$.
Let $A=(a_{ij})$ be an  invertible $k\times k$ integer matrix (hence $\det A= \pm 1$),
and let $v=(v_2,\dots,v_k) \in \R^{k-1}$. Then $ A\cdot v = w =(w_2,\dots,w_k)$
is given by
$$w_j = {
{ a_{1j} + a_{2j}v_2 +\cdots +a_{nj}v_k }
\over
{a_{11} + a_{21}v_2 +\cdots+ a_{n1}v_k}
}, \quad 2 \leq j \leq k.$$
We shall say that the vectors $v,w$ are $GL(k,\Z)$-related.
\end{rema} 
 
\section{Codimension one linear foliations on the torus }
\begin{rema}
{\rm Theorem \ref{TEO} in this section  was obtained long before us by P. ~Iglesias and G.~Lachaud in \cite{IGLESIASLACHAUD}. We are indebted to Professor Iglesias for pointing us the existence of their article, which we were unaware when we wrote our paper.}
\end{rema}

 \subsection{Classification}
Let $\mathcal{F}$ be a codimension one linear foliation with dense leaves on the torus $T^{m+1} =
\R^{m+1} /\Z^{m+1}$.
 It is given by a closed $1$-form
$\omega = dy -\xi_1 dx_1 - \cdots -\xi_n dx_m$, where
at least one of the real numbers   $\xi_1,\dots,\xi_m$ is not rational. The holonomy group
$\Gamma = \left<1,\xi_1,\dots,\xi_m\right>$ is  a dense finitely generated subgroup
of $\R$, hence a free abelian group of rank $k$. 

Our first proposition is a classification of the spaces of leaves, which generalizes the
corresponding result of P. Donato and P. Iglesias for $T^2$ \cite{DONATOIGLESIAS}.

\begin{prop}  \label{ESTAB} Let $\mathcal{F}^\prime$ be another codimension one
linear foliation with dense leaves on $T^{m+1}$, and let $\Gamma^\prime
\subset
\R$ be its holonomy group.  Then the
spaces of leaves
$\R/\Gamma$ and
$\R/\Gamma^\prime$ are diffeomorphic if and only if $\Gamma$ and $\Gamma^\prime$ have equal ranks,  and for any respective basis
$\Gamma=\left\langle\alpha_1, \dots,
\alpha_k\right\rangle$, $\Gamma^\prime=\left\langle\beta_1, \dots, \beta_k\right\rangle$, $k\leq m+1$, 
the vectors $\alpha_1\inv(\alpha_2,\dots,\alpha_k)$ and $\beta_1\inv(\beta_2,\dots,\beta_k)$
are $GL(k,\Z)$-related.
\end{prop} 

In other words, there is some integer matrix $A =(a_{ij}) \in
GL(k,\Z)$ such that
$$
\beta_1\inv\beta_j = ( a_{1j}\alpha_1 + \cdots +a_{kj}\alpha_k)/d, \quad 2\leq j \leq k,
$$
where $d = a_{11}\alpha_1 + \cdots +a_{k1}\alpha_k$.

\proof 
Suppose that  the spaces of leaves $\R/\Gamma$ and $\R/\Gamma^\prime$ are diffeomorphic
as diffeological spaces. By the lifting property  proved in Theorem \ref{LIFT},   there is a diffeomorphism
 $\varphi \colon \R
\to \R$ sending $\Gamma^\prime$ into $\Gamma$. Then $\Gamma$ and $\Gamma^\prime$ have equal ranks.

Moreover, the associated linear automorphism $\theta  =-\varphi(0) +\varphi$ of $\R$ (see Lemma
\ref{MOR}) must be a  homothety
$\theta(t) =\lambda t$ for some $\lambda\in\R$, $\lambda\neq 0$. Since $\theta$
restricted to $\Gamma^\prime$ is an isomorphism into $\Gamma$, we have
$$\lambda \beta_j = a_{1j}\alpha_1 + \cdots + a_{mj}\alpha_k, \quad 1\leq j \leq k,$$
with $a_{ij}\in\Z$. If we take $d=\lambda\beta_1$   the result follows. The converse is immediate.
\qed

\begin{rema}
Notice that if $\Gamma, \Gamma^\prime$ have equal rank $k$, then $\R/\Gamma$ and $\R/\Gamma^\prime$ are
always isomorphic as groups. In fact, since $\R$ is a
$\Q$-vector space, take
$$V=\left<\alpha_1,\dots,\alpha_k,\beta_1,\dots,\beta_k\right>_\Q,$$ which is a
$\Q$-vector space of finite dimension $\geq k$, and write $\R = V\oplus W$ (by Zorn's lemma). Since
$\{\alpha_1,\dots,\alpha_k\}$ and $\{\beta_1,\dots,\beta_k\}$ are linearly independent sets over
$\Q$, there exists $f\colon V \to V$ such that $f(\alpha_i)=\beta_i$, $1\leq i \leq k$, and we can extend
it to $\R$ by putting $f=\id$ on $W$.
\end{rema} 

\begin{coro} [\cite{DONATOIGLESIAS}] \label{DICLAS}
For the torus $T^2$, two irrational linear flows $\mathcal{F}_\alpha$, $\mathcal{F}_\beta$
with holonomy groups $\Gamma_\alpha =\left<1,\alpha\right>$,  $\Gamma_\beta =\left<1,\beta\right>$ have
diffeomorphic (and isomorphic) spaces of leaves if and only if there is some integer matrix
$$A= \left(
\begin{matrix}
a &c \cr
b & d \cr
\end{matrix}
\right)$$
such that
$$\beta = {c + d\alpha \over a + b\alpha}.$$
\end{coro}

\subsection{The group of diffeomorphisms}
\begin{theo}  [\cite{IGLESIASLACHAUD}] \label{TEO} For any codimension one linear foliation on the torus $T^{m+1}$ with
dense holonomy group $\Gamma$, the group $\AUT{R}{\Gamma}$ is isomorphic to
$\Z_2 \times \Z \times
\stackrel{r)}\cdots
\times \Z$, $r\leq m$.
\end{theo}  

A very elegant proof   is given in the paper  \cite{IGLESIASLACHAUD} from Iglesias and Lachaud. A more constructive proof has been done by the third author of the present paper in \cite{SOTELO} as an application of the general formula in subsection \ref{FORMULA}. 
We summarize it.
\proof 
Let $B=\{\alpha_1,\dots,\alpha_k\}$ be a basis of the free abelian group $\Gamma$. We can suppose without loose
of generality that
$\alpha_1=1$. We have to consider three different
cases.
\begin{enumerate}
\item
If the basis $B$ is an algebraic basis over the rationals of some extension $\Q(\beta)$ of finite degree
over $Q$,  let $A$ be the integral closure of $\Q(\beta)$ over $\Z$ \cite{STEWARTTALL}.  Since $\Gamma$ is a $\Z$-module of finite type
with the same rank that $A$, we must have $\Gamma \subset ({1/q})A$ for some positive integer $q$.  It
follows  that
$\AUT{\R}{\Gamma} \subset \AUT{\R}{A} = U(A)$. Hence, by Dirichlet's units theorem, $\AUT{\R}{\Gamma} = \Z_2
\times
\Z
\times
\stackrel{r)}\cdots
\times \Z$, with $0\leq r \leq s+t-1$, where $t$ is the number of  real roots
of the  minimal polynomial of $\beta$ over $\Q$, and  $2s$ is the number
of  complex roots. The
torsion subgroup is
$\Z_2$ because $\pm 1$ are the only real roots of the unit.

\item
If the basis $B$ is not the basis of an  algebraic extension, but all its elements are algebraic,
since $1\in B$, any $\lambda \in \R-\{0\}$ inducing an automorphism of $\Gamma$ must verify
$\lambda^p \in \Gamma$ for all $p\in \Z$. Since the elements of $B$ are algebraic and in a finite number,
we can consider the minimal algebraic extension $\Q(\beta)$ of finite degree over $\Q$ containing $B$.
Then  $\AUT{\R}{\Gamma} \subset U(A)$ as before.
\item
Finally, if the basis $B$ is not the basis of any algebraic extension, and contains transcendent elements,
let $\Gamma_A\subset\Gamma$ be the subgroup of elements of $\Gamma$ which are algebraic numbers over $\Q$. Choose
a basis $B_A$ of $\Gamma_A$ and complete it to a basis $B$ of $\Gamma$. Let $B_T=B-B_A$ the trascendent part of
$B$.
Let
$\lambda\in
\AUT{\R}{\Gamma}$.
Again $0\neq \lambda^p \in \Gamma = \left<B\right>_\Z$ for all $p\in \Z$. Since $\Gamma$ has rank $k$, the
elements
$1,\lambda,\dots,\lambda^k$ must be linearly dependent over $\Z$, hence $\lambda$ is an algebraic
number over $\Q$.  Then $\lambda \in \left<B_A\right>$. Let  $\tau\in B$ be a transcendent element. Then
$\lambda\tau$ is transcendent too (may be with an algebraic part). On the other hand,
 $\lambda$ induces an
automorphism
$\lambda_A \in GL(k,\Z)$ of the free abelian group
$\left<B_A\right>$ of rank $k<m$, because the product of algebraic numbers is   algebraic.
Moreover
$\det \lambda_A = \pm 1$.
That means that
the matrix associated  to $\lambda$ with respect to $\{B_A,B_T\}$ has the form
$$\left(
\begin{matrix}
\lambda_A & \star \cr
0 & \lambda_T \cr
\end{matrix}
\right)
$$
where $\lambda_T \in GL(m-k,\Z)$.  In any case, $\AUT{\R}{\Gamma} \subset \AUT{\R}{\left<B_A\right>}$
and the same arguments that above apply.
\end{enumerate}
\qed

\begin{coro}  Let $\mathcal{F}$ be a   codimension one linear foliation on the torus $T^{m+1}$ with dense leaves. Then
$$\DIFF{T^{m+1}/\mathcal{F}} = (\Z_2 \times \Z \times \stackrel{r)}\cdots
\times
\Z)\ltimes ( T^{m+1}/\mathcal{F}), \quad r\leq m.$$
\end{coro} 

\begin{coro} [\cite{DONATOIGLESIAS}]\label{DONIGL} Let $\mathcal{F}_\alpha$ be an irrational flow on the torus $T^2$.
Then
$\DIFF{T^2/\mathcal{F}_\alpha}$ is isomorphic to
\begin{enumerate}
\item
$(Z_2 \times \Z) \ltimes (T^2/\mathcal{F}_\alpha)$ if $\alpha$ is a quadratic number;
\item
$(Z_2 ) \ltimes (T^2/\mathcal{F}_\alpha)$ otherwise.
\end{enumerate}
\end{coro} 


\subsection{Examples}

\begin{exem} 
{\rm
Let $B = 1 \cup B_T$, with $B_T$ a finite set of transcendent elements (chosen like in the proof of
Theorem \ref{TEO}), such that
$B$ is a linearly independent set over $\Q$. We have seen that any $\lambda\neq 0$ inducing an automorphism
of $\Gamma = \left< B \right>$ must be in $\left< B_A \right> = \Z$ and induce an automorphism of $\left< B_A
\right> = \Z$. Then $\lambda = \pm 1$, so $\AUT{\R}{\Gamma}=\Z_2$. }
\end{exem} 

\begin{exem} {\rm
Let $B=\{ 1, 2^{1/3}, 2^{2/3}\}$. Then it can be proved that  the integral closure of $\Q(2^{2/3})$ over
$\Z$ is exactly $\Gamma = \left< B \right> $. The minimal monic polynomial of $2^{2/3}$ has a real
root  and two
complex roots, that is $t=1,s=1$. Then $\AUT{\R}{\Gamma}= U(\Gamma) = Z_2 \times Z$.
}\end{exem} 

\begin{exem} 
{\rm  Let $B=\{1, \sqrt{2}, e, \sqrt{2}e \}$. The algebraic part of $\Gamma$ is $\Z[\sqrt{2}]$. Then
$B_A=\{1,\sqrt{2}\}$
  can be completed with $B_T=\{e\}$. The unit group of $\Z[\sqrt{2}]$ is $U=\{\pm 1\} \times
\{(1+\sqrt{2})^n\colon n\in \Z\}$. Then if $\lambda \in U$ we have that $\lambda e$ is an automorhism of
$\left< B_T \right>$
if   $\lambda = \pm 1$. Hence
$\AUT{\R}{\Gamma} = \Z_2\times\Z$.
}\end{exem} 

 As in the latter (quadratic) example, the explicit computation of a basis for the diffeomorphism group
would require to solve some algebraic equations (e.g. the Pell-Fermat equation as is done in
\cite{SOTELO}). We shall prove that excepting quadratic foliations,   linear foliations have
no transverse diffeomorphisms other that $\pm \id$.

\section{Linear flows on the torus $T^{n+1}$}\label{FLOWS}
Since a dense flow must have simply connected leaves, by duality we can apply the results of the preceding section to linear {\it flows}.
The importance of linear flows on the torus comes from the following result of P. Caron and Y. Carri\`ere.
\begin{theo}  [\cite{CARON}]  Let $M$ be a compact manifold of dimension $n+1$, endowed with a dense $G$-Lie flow.
Then
$G=\R^n$, $M$ is diffeomorphic to the torus $T^{n+1}$, and the given foliation
is conjugate to a linear one.
\end{theo}

Let $T^{n+1} = \R^{n+1}/\Z^{n+1}$. We take coordinates $(x,y_1,\dots,y_n)$ in $\R\times\R^n$. Let $\omega$
be the closed
$1$-form on the torus, with values on
$\R^n$, given by
$$\omega =(dy_1-\alpha_1dx,\dots,dy_n-\alpha_n dx).$$

\noindent Let $\mathcal{F}_\alpha$ be the corresponding flow. The holonomy group is
$\Gamma =
\left<e_1,\dots,e_n,\alpha\right>
\subset \R^n$, where
$\alpha=(\alpha_1,\dots,\alpha_n)$ and $e_1,\dots,e_n$ is the canonical basis. We shall denote $T_\alpha$ the space of leaves.

As we know,  the abelian group $\Gamma = \left<e_1,\dots,e_n,\alpha\right>$ is dense in $\R^n$ if and only if the
free abelian   group
$\Gamma^\star=\left<1,\alpha_1,\dots,\alpha_n\right>\subset \R$, corresponding to the orthogonal foliation,  has rank $n+1$.
Then the group  $\AUT{\R^n}{\Gamma}$ is isomorphic to $\AUT{\R}{\Gamma^\star}$. In fact, 
  each automorphism in $\AUT{\R^n}{\Gamma}$ corresponds to an integer matrix $A$
 such that $A\cdot \alpha =\alpha$. So we must compute the stabilizer of $\alpha$ for the action
of $GL(n+1,\Z)$. This is exactly what we did in Proposition \ref{ESTAB}, applied to the groups
$\Gamma_\alpha^\star$ and $ \Gamma_\beta^\star$.

\begin{coro} 
The diffeomorphism group $\DIFF{T^{n+1}/\mathcal{F}}$ is isomorphic to $(\Z_2 \times \Z \times \stackrel{r)}{\dots} \times
\Z)
\ltimes (T^{n+1}/\mathcal{F} )$. Moreover $r\leq n$.
\end{coro}

\subsection{Classification}
Again we can classify the spaces of leaves, this time for linear flows on $T^{n+1}$.

\begin{prop}  \label{RELAC} Let $\mathcal{F}_\alpha$, $\mathcal{F}_\beta$ be two dense linear  flows  on the torus $T^{n+1}$,
respectively associated to the vectors $\alpha=(\alpha_1,\dots,\alpha_n)$, $\beta=(\beta_1,\dots,\beta_n)$.
Then their spaces of leaves $T_\alpha$, $T_\beta$ are diffeomorphic if and only if
$\alpha$ and $\beta$ are $GL(n+1,\Z)$ related.
\end{prop} 

\proof  We know that the existence of a diffeomorphism is equivalent to that of some
linear automorphism $\varphi= \left(f_{ij}\right)$ of $\R^n$ such that $\varphi(\Gamma_\alpha) =\Gamma_\beta$.
Since $\varphi(e_j)\in \Gamma_\beta$, we have
\begin{equation}\label{BASE}
(f_{1j},\dots,f_{nj}) = (c_{1j}+b_j\beta_1,\dots,c_{nj}+b_j\beta_n), \quad 1\leq
j\leq n,
\end{equation}
 with $c_{ij}, b_j \in \Z$.
Since $\varphi(\alpha)\in \Gamma_\beta$, we also have
\begin{equation}\label{VECT}
f_{i1}\alpha_1 + \cdots + f_{in}\alpha_n = c_i + b_{n+1}\beta_i, \quad 1\leq i \leq n,
\end{equation}
with $c_i,  b_{n+1} \in \Z$.
From equations (\ref{BASE}) and (\ref{VECT}) we obtain
$$
\beta_i = { {c_i- c_{i1}\alpha_1- \cdots -c_{in}\alpha_n} \over {-b_{n+1}+b_1\alpha_1 + \cdots +
b_n\alpha_n}}, \quad 1\leq i \leq n,
$$

that is $\beta = A\cdot \alpha$ for
$$A= \left(
\begin{matrix}
-b_{n+1} & c_1 & \cdots & c_n \cr
b_1 & -c_{11} &\cdots & -c_{n1}&\cr
\vdots & \vdots &&\vdots\cr
b_n &  -c_{1n} &\cdots &  -c_{nn}\cr
\end{matrix}
\right) \in GL(n+1,\Z).
$$
\qed

When $n=1$ the latter result was proven in \cite{DONATOIGLESIAS} (see also Corollary \ref{DICLAS}).
\section{Transcendent and quadratic foliations}
Our next generalization of the results of \cite{DONATOIGLESIAS} involves arbitrary linear foliations of dimension
$m\geq 1$ and codimension $n\geq 1$ on the torus $T^{m+n}$. Such a foliation is determined by  some
linear subspace
$V\subset \R^{m+n}$ of dimension $m$. By taking dual coordinates, we have a closed $1$-form $\omega
=(\omega_1,\dots,\omega_n)$ with values in $\R^n$ which defines the foliation.

\subsection{Transcendent foliations}
Let us say that the foliation is {\it transcendent} if the invariant subspace $V$ admits a basis
\begin{eqnarray*}
v_1 &=& (\alpha_1^1,\dots,\alpha_m^1,\beta_1^1,\dots,\beta_n^1)\\
&\vdots & \\
v_m &= &(\alpha_1^m,\dots,\alpha_m^m,\beta_1^m,\dots,\beta_n^m)
\end{eqnarray*}
\noindent such that   all the coordinates $\alpha_i^j, \beta_i^j$ are
algebraically independent over
$\Q$.  The basis
$\{v_1,\dots,v_m\}$ will be called a {\it transcendence basis}. Notice that a transcendent foliation has
 dense leaves, because,  accordingly to what we have seen in the preceding sections,  the flow generated by any of the $v_i$ is
dense.

Transcendent foliations appeared in B.~Herrera's thesis \cite{HERRERA}. Our definition
is easily seen to be equivalent to the original one (where the last coordinates of the $v_i$'s are
supposed to be $1$).

\begin{theo}  [\cite{HERRERA}]
The only transverse diffeomorphisms of a transcendent foliation
are $\pm \id$.
\end{theo}  

 We shall generalize this result by dualizing the definition above.

\begin{lemm} 
The subspace $V\subset \R^{m+n}$ generates a transcendent foliation if and only if its
orthogonal subspace $W=V^\perp$ generates a transcendent foliation.
\end{lemm} 

\proof If $n=1$,  let us take a transcendence  basis   of $V$, and the
vector
$$
w =  v_1\wedge \cdots \wedge v_m= \det \left(
\begin{matrix}
 e_1 & \dots &e_m & e_{m+1} \cr
\alpha_1^1 & \dots &   \alpha_m^1 & \beta_1^1 \cr
\vdots & &\vdots &\vdots \cr
\alpha_1^m & \dots &  \alpha_m^m & \beta_1^m \cr
\end{matrix}
\right),
$$
which is a basis of $W$,
and whose  coordinates remain algebraically independent over $\Q$.

When $V$ has codimension $n\geq 2$, we can complete any transcendence basis $\{v_1,\dots,v_m\}$ of $V$
to a  transcendence
basis $\{v_1,\dots,v_m, u_1,\dots,u_n \}$ of $\R^{m+n}$, because the set of real
numbers which are algebraically dependent over any extension $\Q(t_1,\dots,t_r)$, $t_1,\dots,t_r \in \R$,
is a countable set.  Let us consider the vectors
$$w_j = v_1 \wedge \cdots \wedge v_m \wedge u_1 \wedge \cdots \wedge \widehat{u_j} \wedge \cdots
\wedge u_n, \quad 1\leq j \leq n.$$
Then $\{w_1,\dots,w_n\}$ is a transcendence basis of $W$.\qed

\begin{coro}  The foliation $\mathcal{F}$ is transcendent if and only if the coordinates of some closed form
$\omega=(\omega_1,\dots,\omega_n)$ defining it are a set of algebraically independent numbers over $\Q$.
\end{coro} 

\proof The coordinates of $\omega$ with respect to the dual basis
$$dx_1,\dots,dx_m,dy_1,\dots,dy_n$$ just correspond to
the orthogonal subspace
$W$. \qed

\begin{prop}  The foliation is transcendent if and only if it can be defined by $1$-forms
$$
\omega_j = dy_j +\beta_j^1 dx_1 + \cdots +\beta_j^m dx_m, \quad 1\leq j \leq n,
$$
such that the coefficients $\{\beta_i^j\}$, $1\leq i\leq m$, $1\leq j \leq n$, are algebraically
independent.
\end {prop}

\proof Since   the foliation is transcendent, we have $\det (\alpha_i^j) \neq 0$ for the
coefficients of the transverse part of the forms defining the foliation. The corresponding change
of basis respects the algebraic independence of the tangential coefficients. \qed

\subsection{Quadratic foliations}
The considerations of the preceding subsection suggest the following definition.
\begin{defi} 
Let us say that  a codimension $n$ linear foliation  $\mathcal{F}$  on the torus $T^{m+n}$ is  {\em non-quadratic}
whenever it can be defined by a closed
$\R^n$-valued $1$-form
$\omega =(\omega_1,\dots,\omega_n)$, where
$$\omega_j = dy_j + \beta_j^1dx_1 + \cdots +\beta_j^m dx_m, \quad 1\leq j \leq n,$$
such that any polynomial with $mn$ variables and rational coefficients that annihilates the
family of coefficients $\{\beta_j^i\}$ has its degree greater than $2$.
\end{defi} 

Obviously a transcendent foliation is not quadratic.

\begin{exem} {\rm The foliation on $T^3$ defined by the $2$-subspace $V\subset \R^3$ with basis
$v_1=(1,0,2^{1/3})$, $v_2=(0,1 , 3^{1/3})$  is neither transcendent nor quadratic.
}\end{exem} 


\begin{theo}   \label{RIGID}Let $\mathcal{F}$ be a non-quadratic linear foliation.    Then $\AUT{\R^n}{\Gamma} =
\{\pm \id\}$ and $\DIFF{T^{m+n}/\mathcal{F}} = \Z_2\ltimes
(T^{m+n}/\mathcal{F}).$
\end{theo}  

\proof The holonomy group is $\Gamma=\left<e_1,\dots,e_n, \beta_1,\dots,\beta_m,\right> \subset \R^n$, where
$\beta_k = (\beta_1^k,\dots,\beta_n^k)$. In the same way as we did in Proposition \ref{RELAC},
one can prove that the matrix $\left(f_{ij}\right)$ associated to a linear automorphism $\varphi\in
\AUT{\R^n}{\Gamma}$ must verify
$$f_{ij} = c_{ij} + b_1^j\beta_i^1 + \cdots + b_m^j\beta_i^m, \quad 1\leq i,j \leq n,$$
with $c_{ij}, b_i^j \in \Z$ because $\varphi(e_j)\in\Gamma$, $1\leq j\leq n $. From the conditions
$\varphi(\beta_k) \in\Gamma$, $1\leq k\leq m$, we  deduce relations  of algebraic dependence over $\Q$ for
the family $\{\beta_k^j\}$, that is
$$ f_{i1}\beta^k_1 + \cdots + f_{in}\beta^k_n \in \left< 1, \beta_i^j\right>_\Z, \quad 1\leq i \leq k, 1\leq k
\leq m.$$
These relations  are given by a rational polynomial on
$mn$ variables,   having degree not greater than $2$, excepting when
$c_i^j=0$,
$1\leq i \leq m$, $1\leq j \leq n,$ and $m_{ij}=0$ when $i\neq j$.  Hence
the matrix is diagonal with the same number $m\in\Z$ in all entries. Since $\varphi$ induces
an automorphism of $\Gamma$, we deduce analogous conditions for the inverse matrix, so $m=\pm 1$ and
$\varphi =\pm \id$. \qed

For $n=1$ this is again Donato-Iglesias' result cited in Corollary \ref{DONIGL}. For transcendent foliations
it has been proved by B. Herrera in \cite{HERRERA}.

\subsection{Classification}
The proof of the
following theorem is analogous to that of Proposition
\ref{RELAC} and we shall omit it.
\begin{theo}  
Let $\mathcal{F},\mathcal{F}^\prime$ be two non-quadratic linear foliations on $T^{m+n}$ associated
to the coefficients $\beta_1,\dots,\beta_m\in\R^n$ and $\gamma_1,\dots,\gamma_m\in\R^n$ respectively. Then the spaces of
leaves
$T^{m+n}/\mathcal{F}$ and $T^{m+n}/\mathcal{F}^\prime$ are diffeomorphic (in the diffeological sense) if and
only if
$$\gamma = (A+\beta B)\inv (C+\beta D)$$
for integer matrices $A\in \mathcal{M}_{n\times n}(\Z)$, $B\in\mathcal{M}_{m\times m}(\Z)$,
$C\in\mathcal{M}_{n\times m}(\Z)$ and $D\in \mathcal{M}_{m\times n}(\Z)$, where $\beta,\gamma$ are the matrices
whose columns are the given coefficients of the foliation.
\end{theo}  
This could be stated in terms of the action of $GL(m+n,\Z)$ on the grasmannian $G_n^{m+n}$. For $m=1, n=1$
we reobtain Corollary
\ref{DICLAS}.

\nocite{*}
\bibliographystyle{plain}
\bibliography{bibliocopia}
\end{document}